\documentclass[12pt]{article}
\begin{document}
\newcommand{\qed}{\ \
\mbox{\rule{8pt}{8pt}}\vspace{0.3cm}\newline}
\newcommand{\ia}{{\bf I}_{A_i}}
\newcommand{\ba}{\widetilde{{\bf KC}}_i}
\newcommand{\bba}{{\bf KC}_i}
\newcommand{\ra}{\longrightarrow}
\newcommand{\pe}{{\cal P}}
\newcommand{\der}{{\cal DP}_d}
\newcommand{\adp}{{\cal DP}_d^{af}}
\newcommand{\ot}{\otimes}
\newcommand{\rec}{\raisebox{-1ex}{\ $\stackrel{\textstyle{\stackrel{\textstyle{\longleftarrow}}{\longrightarrow}}}{\longleftarrow}$\ }}
\title{Poincar\'e duality for Ext--groups between strict polynomial functors}
\author{Marcin Cha\l upnik
\thanks{The author was supported by the  grant (NCN) 2011/01/B/ST1/06184.}\\
\normalsize{Institute of Mathematics, University of Warsaw,}\\
\normalsize{ul.~Banacha 2, 02--097 Warsaw, Poland;}\\
\normalsize{e--mail: {\tt mchal@mimuw.edu.pl}}}
\date{\mbox{}}
\newtheorem{prop}{Proposition}[section]
\newtheorem{cor}[prop]{Corollary}
\newtheorem{theo}[prop]{Theorem}
\newtheorem{lem}[prop]{Lemma}
\newtheorem{defi}[prop]{Definition}
\newtheorem{defipro}[prop]{Definition/Proposition}
\newtheorem{fact}[prop]{Fact}
\newtheorem{exa}[prop]{Example}
\newcommand{\ka}{{\mbox {\bf k}}}
\newcommand{\ca}{\mbox{$\cal{A}$}}
\maketitle
\begin{abstract}
We study relation between left and right adjoint functors to the precomposition functor. As a consequence we obtain various dualities in the Ext--groups in the category of strict polynomial functors.\\\mbox{}\vspace{0.1cm}
\end{abstract}
\section{Introduction}
Poincar\'e duality patterns in the category ${\cal P}_d$ of strict polynomial functors were observed empirically in many examples [FFSS, C1, C2, C3].
Typically, they occur in Ext--groups involving Frobenius
twisted functors. In the present paper we provide a conceptual explanation of these phenomena.
The main ingredient is an observation (Theorem~2.1) that the left and right adjoint functors to the precomposition with the Frobenius twist are closely related.
Then by applying the contravariant involution $(-)^{\#}$ on ${\cal P}_d$ we obtain various  Poincar\'e--like formulae (Cor.~3.3, Ex.~3.4, Ex.~3.5) for the Ext--groups. Our construction formally
resembles  the construction of Poincar\'e duality for sheaves on algebraic varieties via the Verdier duality. We discuss these analogies in Section 4.\newline
{\em Acknowledgements.} I am grateful to Stanis\l aw Betley whose ongoing work on cubical construction inspired the present article and to Andrzej Weber
for many discussions on sheaves and dualities.
\section{Cohomological shift}
Let ${\cal P}_d$ be the category of strict polynomial functors homogeneous of degree $d$ (see [FS]) over a field {\bf k} of characteristic $p>0$ and
let ${\cal DP}_d$ be the bounded derived category of cohomological complexes over ${\cal P}_d$. We recall from [C4] that the functor
\[{\bf C}:{\cal DP}_d\ra{\cal DP}_{p^id}\]
of taking precomposition with the $i$th Frobenius twist
has right and left adjoint functors
\[{\bf K}^r,{\bf K}^l:{\cal DP}_{p^i d}\longrightarrow {\cal DP}_d.\]
In order to explicitly describe ${\bf K}^r$ and ${\bf K}^l$
we introduce the following family of strict polynomial functors: for
$U\in \mbox{Vect}^{fin}_{{\bf k}}$ we define $\Gamma^{d(i)}_{U^*}\in {\cal P}_{p^id}$ by the formula $V\mapsto \Gamma^{d}(U^*\otimes V^{(i)})$
where $\Gamma^d$ is the divided power functor, $(-)^*$ stands for the {\bf}--linear dual and $(-)^{(i)}$ is the $i$th Frobenius twist. Now we put
\[{\bf K}^r(F)(V):=\mbox{RHom}_{{\cal D}{\cal P}_{p^i d}}(\Gamma^{d(i)}_{V^*},F),\]
and
\[{\bf K}^l(F):={\bf K}^r(F^{\#})^{\#}\]
where $F^{\#}(V):=(F(V^*))^*$ is the Kuhn dual of $F$.
The source of duality phenomena in ${\cal DP}_d$ is the following relation between ${\bf K}^r$ and ${\bf K}^l$.
\begin{theo}
There is a natural isomorphism of functors
\[{\bf K}^l(F)\simeq {\bf K}^r(F)[2(p^i-1)d].\]
\end{theo}
{\bf Proof: } For $U\in \mbox{Vect}^{fin}_{{\bf k}}$ we define $S^{p^id}_U\in {\cal P}_{p^id}$ as $V\mapsto S^{p^id}(U\otimes V)$.
It is well known (see e.g. [FS, Th. 2.10]) that  $\{S^{p^id}_U\}$ is a family of injective cogenerates ${\cal P}_{p^id}$. Therefore in order to prove Theorem 2.1 it suffices to construct a natural in $U$ isomorphism
\[{\bf K}^l(S^{p^id}_U)\simeq {\bf K}^r(S^{p^id}_U)[2(p^i-1)d].\]
This isomorphism  essentially follows from the following computation of Hom/Ext--groups.
\begin{lem}
We have natural in $V,U$ isomorphisms
\begin{enumerate}
\item
$\mbox{Hom}_{{\cal P}_{p^id}}(\Gamma^{d(i)}_{V^*},S^{p^id}_U)\simeq S^d_{U^{(i)}}(V),$
\item
$\mbox{Ext}^*_{{\cal P}_{p^id}}(\Gamma^{d(i)}_{V^*},\Gamma^{p^id}_U)\simeq \Gamma^d_{U^{(i)}}(V)[-2(p^i-1)d].$
\end{enumerate}
\end{lem}
{\bf Proof of the Lemma: }
The first formula immediately follows from the Yoneda lemma [FS, Th. 2.10], so we turn to the proof of the second formula.
We start with the case  $d=1$. We clearly have
\[\mbox{Ext}^*_{{\cal P}_{p^i}}(\Gamma^{1(i)}_{V^*},\Gamma^{p^i}_U)\simeq V\ot \mbox{Ext}^*_{{\cal P}_{p^i}}(I^{(i)},\Gamma^{p^i}_U),\]
hence our task is reduced to showing the following parameterized version of [FFSS,~Th.5.4]:
\[\mbox{Ext}^*_{{\cal P}_{p^i}}(I^{(i)},\Gamma^{p^i}_U)\simeq U^{(i)}.\]
To see this we consider the dual Koszul complex with parameter $U$ (see e.g. [FS, Sect.~4])
\[0\ra \Gamma^{p^i}_U\ra \Gamma^{p^i-1}_U\ot \Lambda^1_U\ra\ldots\ra\Lambda^{p^i}_U\ra 0.\]
Since this complex is exact and $\mbox{Ext}^*_{{\cal P}_{p^i}}(I^{(i)},\Gamma^j_U\ot\Lambda^{p^i-j}_U)=0$ for $0<j<p^i$ [FS, Th.~2.13], we get a natural identification
\[ \mbox{Ext}^*_{{\cal P}_{p^i}}(I^{(i)},\Gamma^{p^i}_U)\simeq \mbox{Ext}^*_{{\cal P}_{p^i}}(I^{(i)},\Lambda^{p^i}_U)[-(p^i-1)].\]
Now we take undualized Koszul complex with parameter $U$:
\[0\ra \Lambda^{p^i}_U\ra \Lambda^{p^i-1}_U\ot S^1_U\ra\ldots\ra S^{p^i}_U\ra 0,\]
and analogously conclude that
\[ \mbox{Ext}^*_{{\cal P}_{p^i}}(I^{(i)},\Lambda^{p^i}_U)\simeq \mbox{Ext}^*_{{\cal P}_{p^i}}(I^{(i)},S^{p^i}_U)[-(p^i-1)]\simeq
U^{(i)}[-(p^i-1)].\]
This gives our formula for $d=1$.\newline

In order to get this formula for $d>1$ we first
consider the groups $\mbox{Ext}^*_{{\cal P}_{p^id}}(\Gamma^{d(i)}_{V^*},(\Gamma^{p^i}_U)^{\otimes d})$.
Then by the K\"unneth formula [FFSS, p.~672] we get
\[\mbox{Ext}^*_{{\cal P}_{p^id}}(\Gamma^{d(i)}_{V^*},(\Gamma^{p^i}_U)^{\otimes d})\simeq
(\mbox{Ext}^*_{{\cal P}_{p^i}}(\Gamma^{1(i)}_{V^*},(\Gamma^{p^i}_U)))^{\otimes d}\simeq\]
\[(V\otimes U^{(1)}[-2(p^i-1)])^{\otimes d}.\]
Then by [FFSS, Th.~5.4] the inclusion $\Gamma^{p^id}_U\longrightarrow (\Gamma^{p^i}_U)^{\otimes d}$ induces a monomorphism on Ext--groups
\[\mbox{Ext}^*_{{\cal P}_{p^id}}(\Gamma^{d(i)}_{V^*},\Gamma^{p^id}_U)\longrightarrow
(\mbox{Ext}^*_{{\cal P}_{p^id}}(\Gamma^{d(i)}_{V^*},(\Gamma^{p^i}_U)^{\otimes d}))^{\Sigma_d}.\]
Since, by [FFSS, Th.~5.4] again, the graded dimensions of both sides are equal, we get the isomorphisms
\[\mbox{Ext}^*_{{\cal P}_{p^id}}(\Gamma^{d(i)}_{V^*},\Gamma^{p^id}_U)\simeq
(\mbox{Ext}^*_{{\cal P}_{p^id}}(\Gamma^{d(i)}_{V^*},(\Gamma^{p^i}_U)^{\otimes d}))^{\Sigma_d}\simeq\]
\[((V\otimes U^{(1)}[-2(p^i-1)])^{\otimes d})^{\Sigma_d}\simeq \Gamma^d_{U^{(i)}}(V)[-2(p^i-1)d].\]
\qed
Now by the injectivity of $S^{p^id}_U$ and the first part of the lemma we get ${\bf K}^r(S^{p^id}_U)\simeq S^d_{U^{(i)}}$. On the other hand, since
$(S^{p^id}_U)^{\#}\simeq \Gamma^{p^id}_{U^*}$, we get
\[H^*({\bf K}^l(S^{p^id}_U))=H^*({\bf K}^r(\Gamma^{p^id}_{U^*}))^{\#}\simeq \Gamma^d_{U^{*(i)}}[-2(p^i-1)d]^{\#}\simeq S^d_{U^{(i)}}[2(p^i-1)d].\]
It remains to construct a natural in $U$ quasi--isomorphism $H^*({\bf K}^l(S^{p^id}_U))\simeq {\bf K}^l(S^{p^id}_U)$. To this end we observe that
$\Gamma^{d(i)}_{V^*}$ has, as an object of ${\cal P}_{p^id}$, homological dimension $2(p^i-1)d$ (e.g. we can take the (non--natural in $V$) Troesch
resolution [Tr]). Thus $H^j({\bf K}^r(F))$ for $F\in{\cal P}_{p^id}$ can be non--trivial only for $j\leq 2(p^i-1)d$. Thus, by a general homological algebra argument,
$\Gamma^{d(i)}_{V^*}$ has a resolution of length $2(p^i-1)d$ in ${\cal P}^d_{p^id}$ acyclic with respect to ${\bf K}^r$. Then when we use such a resolution for computing ${\bf K}^r(\Gamma^{p^id}_{U^*})$ we see that the
${\bf K}^l(S^{p^id}_U)$ is a complex spanning between degrees $-2(p^i-1)d$ and $0$. Hence its cohomology  is concentrated in the lowest degree of the complex.
Thus the quasi--isomorphism $H^*({\bf K}^l(S^d_U))\simeq {\bf K}^l(S^d_U)$ can be realized just as the natural embedding of the lowest cohomology group of a bounded complex.\qed
\section{Duality in Ext--groups}
The applications to Ext--groups follow from an elementary general lemma.
Let us call a \ka--linear abelian category ${\cal A}$ \ka--finite when every object has a finite composition series and Hom--spaces are finite dimensional.
\begin{lem}
Let ${\cal A}$ be a \ka--finite category and let $D:{\cal A}\ra {\cal A}^{op}$ be an equivalence such that all simple objects
are self--dual. Let $P\in {\cal A}$ be projective and self--dual (hence also injective) and let $R:=(\mbox{End}_{{\cal A}}(P))^{op}$.
Assume that $R$ is a symmetric Frobenius \ka--algebra (e.g. it is always the case when $P$ is simple)
Then we  have an isomorphism
\[\mbox{Hom}_{{\cal A}}(P,-)\simeq\mbox{Hom}_{{\cal A}}(-,P)^*\]
of functors from ${\cal A}$ to $R$--mod.
\end{lem}
{\bf Proof: }
Let us denote by $(-)^{\diamond}$ the $R$--dual i.e. for a right $R$--module $M$, $M^{\diamond}:=\mbox{Hom}_R(M,R)$. Then since
$R$ is a symmetric Frobenius \ka--algebra,
$R\simeq R^*$ as  $R$--$R$ bimodules. Hence for any right $R$--module $M$ we have natural isomorphisms  of left $R$--modules:
\[M^*\simeq\mbox{Hom}_R(M,R^*)\simeq\mbox{Hom}_R(M,R)=M^{\diamond}.\]
This shows that the functors $\mbox{Hom}_{{\cal A}}(-,P)^*$ and $\mbox{Hom}_{{\cal A}}(-,P)^{\diamond}$ are isomorphic. Thus it suffices to construct
an isomorphism
\[\mbox{Hom}_{{\cal A}}(P,-)\simeq\mbox{Hom}_{{\cal A}}(-,P)^{\diamond}.\]
We consider  a natural in $X$,  $R\ot R^{op}$ equivariant map
\[\mu_X:\mbox{Hom}_{{\cal A}}(X,P) \otimes_R\mbox{Hom}_{{\cal A}}(P,X)\ra R. \]
given by the  composition: $\mu_X(\phi\ot \psi):=\phi\circ \psi$.
This allows us to define  transformation
\[\Psi: \mbox{Hom}_{{\cal A}}(P,-)\ra\mbox{Hom}_{{\cal A}}(-,P)^{\diamond}\]
by the formula
\[\Psi_X(\psi):=\mu_X(-\ot\psi).\]
Since objects of ${\cal A}$ have finite composition series and both functors are exact, it suffices to show that $\Psi_S$ is an isomorphism for any simple $S\in{\cal A}$.
Now, since $S$ and $P$ are self--dual, $\mbox{Hom}_{{\cal A}}(P,S)$ and $\mbox{Hom}_{{\cal A}}(S,P)^*\simeq
\mbox{Hom}_{{\cal A}}(S,P)^{\diamond}$ have (finite and) equal dimension over \ka. Therefore it suffices
to show that $\Psi_S$ is a monomorphism which is equivalent to showing that $\mu_S$ is right non--degenerate.
To this end  we observe that, since $S$ is simple, any non--trivial $\psi: P\ra S$ is epimorphic and any non--trivial $\phi:S\ra P$ is monomorphic.
Hence their composition is non--trivial. \qed
{\bf Remark:} It is not obvious whether the assumption on $R$ is really necessary. In fact, it is easy to show that
$\mbox{Hom}_{{\cal A}}(P,X)$ and $\mbox{Hom}_{{\cal A}}(X,P)^*$ have equal $\ka$--dimension without any assumptions on $R$.
However, since this assumption is satisfied for all examples we have in mind, we leave this question unanswered.\newline

Lemma 3.1 together with a shift phenomenon discussed in the previous section produce the Poincar\'e duality in our context.
\begin{theo}
Let $P\in{\cal P}_d$ be either simple, projective  or $P=I^d$ (the d--th tensor power functor).
Then for any $s\geq 0$ we have a natural in $F$ isomorphism
\[ \mbox{Ext}^s_{{\cal P}_{p^id}}(P^{(i)}, F)\simeq \mbox{Ext}^{2(p^i-1)d-s}_{{\cal P}_{p^id}}(P^{(i)}, F^{\#})^*.\]
\end{theo}
{\bf Proof:} We shall apply Lemma 3.1 to the category ${\cal P}_d$ with the Kuhn duality $(-)^{\#}$ as $D$. Indeed, it is well known that all simples are self--dual
and that $I^d$ is self--dual and projective.
Also we obtain by standard  computations (see e.g. [FF, Th.~1.8]) that $\mbox{Hom}_{{\cal P}_d}(I^d,I^d)\simeq\ka[\Sigma_d]$ which as a group algebra is symmetric Frobenius.
 Then by using $\{{\bf C},{\bf K}^r\}$ adjunction and projectivity of $P$ we get
\[\mbox{Ext}^s_{{\cal P}_{p^id}}(P^{(i)}, F)\simeq \mbox{HExt}^s_{{\cal P}_{p^id}}(P,{\bf K}^r(F))\simeq H^s(\mbox{Hom}_{{\cal P}_d}(P,{\bf K}^r(F))).\]
Similarly, using this time the injectivity and self--duality of $P$ and additionally Theorem 2.1  we obtain
\[\mbox{Ext}^{2(p^i-1)d-s}_{{\cal P}_{p^id}}(P^{(i)}, F^{\#})^*\simeq \mbox{Ext}^{2(p^i-1)d-s}_{{\cal P}_{p^id}}(F,P^{(i)})^*\simeq
\mbox{HExt}^{2(p^i-1)d-s}_{{\cal P}_{d}}({\bf K}^l(F),P)^*\simeq\]\[ H^{2(p^i-1)d-s}(\mbox{Hom}_{{\cal P}_d}({\bf K}^l(F),P))^*\simeq H^{-s}(\mbox{Hom}_{{\cal P}_d}({\bf K}^r(F),P))^*\simeq H^{s}(\mbox{Hom}_{{\cal P}_d}({\bf K}^r(F),P)^*).\]
Now our claim follows from Lemma 3.1.\qed
In order to obtain concrete examples we recall some basic facts concerning simple objects in ${\cal P}_d$ (see e.g. [Ma]).
The family of simples in ${\cal P}_d$
is indexed by the set of Young diagrams of weight d.
Since ${\cal P}_d$ is a highest weight category, if certain block consists of a single simple then this simple is projective.
Such simples are labeled by Young diagrams which have $p$--weighting $e$ and are  $p^{e+1}$--cores for some $e\geq 0$ [Ma, Sect.~5]. Therefore we get
\begin{cor}
Let $F_{\mu}$ be a simple functor associated to a Young digram $\mu$ of weight d
 which has $p$--weighting $e$ and is a  $p^{e+1}$--core
and let $F_{\lambda}$ be a simple functor associated to a Young diagram
$\lambda$ of weight $p^id$. Then for any $s\geq 0$
\[ \mbox{Ext}^s_{{\cal P}_{p^id}}(F_{\mu}^{(i)}, F_{\lambda})\simeq \mbox{Ext}^{2(p^i-1)d-s}_{{\cal P}_{p^id}}(F_{\mu}^{(i)}, F_{\lambda})^*.\]
\end{cor}
The importance of these Ext--groups was indicated already in [C3, Sect.~5; C4, Sect.~4].
They seem to  play a fundamental role in the structure of ${\cal DP}_{p^id}$, at least they are building blocks for many other Ext--groups.\newline
Let us point out the simplest instance of Corollary 3.3 (or the case of $d=1$ for $P=I^d$ in Theorem 3.2).
\begin{exa}
For any Young diagram $\lambda$ of weight $p^i$ and any $s\geq 0$
\[ \mbox{Ext}^s_{{\cal P}_{p^i}}(I^{(i)}, F_{\lambda})\simeq \mbox{Ext}^{2(p^i-1)-s}_{{\cal P}_{p^i}}(I^{(i)}, F_{\lambda})^*.\]
\end{exa}
We can obtain further Poincar\'e dualities between Ext--groups by combining  Theorem 3.2 with the Koszul duality [C2]. We recall that the Koszul duality functor $\Theta: {\cal DP}_d\ra {\cal DP}_d$ is a self--equivalence which takes the Schur functor $S_{\lambda}$ to the Weyl functor $W_{\widetilde{\lambda}}$ associated to the conjugate Young diagram $\widetilde{\lambda}$ and takes $I^{d(i)}$ to $I^{d(i)}[-p^id]$. Taking into account the fact that $W_{\widetilde{\lambda}}^{\#}=S_{\widetilde{\lambda}}$ we get
\begin{exa}
For any Young diagram $\lambda$ of weight $p^id$ and any $s\geq 0$
\[ \mbox{Ext}^s_{{\cal P}_{p^id}}(I^{d(i)}, S_{\lambda})\simeq \mbox{Ext}^{(p^i-1)d-s}_{{\cal P}_{p^id}}(I^{d(i)}, S_{\widetilde{\lambda}})^*.\]
\end{exa}
This duality for $i=1$ can be nicely interpreted combinatorially. Namely, it was shown in [C3, Sect.~5] that $\mbox{Ext}^*_{{\cal P}_{pd}}(I^{d(1)}, S_{\lambda})$
has a basis labeled by the set of fillings of $\lambda$ by rim p--hooks (called there ``slicings''). There is of course a bijection between  slicings
of $\lambda$ and of $\widetilde{\lambda}$ and when we look carefully how the Ext--degree depends on the shapes of rim p--hooks in the slicing we recover
our Example 3.5.

\section{Geometric analogy}
In this section we briefly discuss formal similarities between our construction and that of Poincar\'e duality for sheaves.
The categorical essence  of the sheaf construction (see e.g [Ha]) can be extracted
 as follows. We have a triangulated \ka--linear category ${\cal T}$ in which we would like to have duality and another (usually simpler) triangulated \ka--linear category ${\cal U}$ which will be used for producing cohomological functors. We have a pair of adjoint functors $\{i^*,i_*\}$ between them i.e. $i_*:{\cal T}\ra{\cal U}$,  $i^*:{\cal U}\ra{\cal T}$. Both categories are equipped with contravariant involutions (in both categories denoted by $D$) and $i_*$ commutes with $D$. Then, automatically, the functor $i^!:=D\circ i^*\circ D$ is right adjoint to $i_*$. At last we have a self--dual $P\in{\cal U}$ such that  the functors $(P,-)_{{\cal U}}$ and $(-,P)_{{\cal U}}^*$ are isomorphic.
Then we define a cohomolocial functor ${\bf H}:{\cal T}\ra \ka\mbox{--mod}_{gr}$ defined by the formula ${\bf H}(X):=(P,i_*(X))_{{\cal U}}$ and a homologicaal
functor ${\bf H}':{\cal T}\ra \ka\mbox{--mod}_{gr}^{op}$ defined by  ${\bf H}'(X):=(i_*(X),P)_{{\cal U}}$. By our assumption, ${\bf H}\simeq {\bf H}'^*$. Now we have
\[{\bf H}'(X)=(i_*(X),P)_{{\cal U}}\simeq (i_*(X),D(P))_{{\cal U}}\simeq (X,i^!(D(P)))_{{\cal T}}\simeq\]\[ (X,D(i^*(P)))_{{\cal T}}\simeq (i^*(P), D(X))_{{\cal T}}
\simeq (P,i_*(D(X))_{{\cal U}}={\bf H}(D(X)).\]
Hence we have a natural isomorphism
\[{\bf H}(X)\simeq {\bf H}'(D(X))\simeq {\bf H}(D(X))^*.\]
If additionally there is an isomorphism $i^!\simeq i^*[n]$ for some $n>0$ then for any self--dual $Y\in{\cal U}$ (e.g. $Y=P$) we have
\[D(i^*(Y))\simeq i^!(D(Y))\simeq i^!(Y)\simeq i^*(Y)[n].\]
Thus we get
\[{\bf H}(i^*(Y))\simeq {\bf H}(i^*(Y))^*[-n].\]
As we have mentioned, an example motivating this general construction comes from sheaf theory. Namely, let $M$ be a smooth projective variety over complex
numbers. We put ${\cal T}$ to be the derived category of the category of constructible sheaves of {\bf C}--vector spaces on $M$ and ${\cal U}$ to be the analogous category on a one--point space. Then we take $P$ to be {\bf C} put in degree 0, $i_*$ and $i^*$ respectively the (derived) direct and inverse image functors for
the map $M\ra *$. At last, $D$ is the Verdier duality functor. These data satisfy all the conditions we formulated in the  abstract context, therefore we can apply the machinery and obtain duality. Since $i^*({\bf C})$ is the constant sheaf on $M$, we recover the classical Poincar\'e duality.\newline
Thus we see that the situation we consider in the present paper is  dual to that described above. Namely, we have two triangulated categories
${\cal T}:={\cal DP}_{p^id}$,  ${\cal U}:={\cal DP}_{d}$ equipped with dualities $D:=(-)^{\#}$ but the functors go into another directions:
$i_*:={\bf C}:{\cal U}\ra {\cal T}$, $i^*:={\bf K}^l:{\cal T}\ra {\cal U}$. Let us look at this dual situation in abstract terms. This time
we have a cohomological functor
${\bf H}:{\cal T}\ra \ka\mbox{--mod}_{gr}$ defined by the formula ${\bf H}(X):=(P,i^!(X))_{{\cal U}}$ and a homological one
defined by ${\bf H}'(X):=(i^!(X),P)_{{\cal U}}$ and by the assumption on $P$, ${\bf H}\simeq {\bf H}'^*$.
By using adjunctions we get
\[{\bf H}'(X):=(i^!(X),P)_{{\cal U}}\simeq (i^!(X),D(P))_{{\cal U}}\simeq(P,D(i^!(X)))_{{\cal U}}\simeq(P,i^*(D((X))))_{{\cal U}}.\]
Here, already at this point we need an additional assumption $i^!\simeq i^*[-n]$. With it we get
\[(P,i^*(D((X))))_{{\cal U}}\simeq(P,i^!(D((X))))_{{\cal U}}[n]\simeq {\bf H}(D(X))[n].\]
Hence this time we obtain the following form of duality
\[{\bf H}(X)\simeq {\bf H}'(X)^*\simeq {\bf H}(D(X))^*[-n].\]
When we apply it to $X=i_*(Y)$ for a self--dual $Y\in{\cal U}$, then, since $D(i_*(Y))\simeq i_*(Y)$, we get
\[{\bf H}(i_*(Y))\simeq {\bf H}^*(i_*(Y))[-n].\]
Thus, eventually, we have arrived at a similar formula to that  in the dual setting, but the  point is that (as we have seen in Corollary 3.3) we have much more interesting self--dual objects in ${\cal T}$ than only those coming from ${\cal U}$. Thus in some respect this approach provides  more easily accessible examples than the dual one.\newline
On the other hand we can formally turn our setting into  the previous one. Namely, we can put $i_*:={\bf K}^r$ and $i^*:={\bf C}$. Also it is easy
to see that $i^!:={\bf C}[2(p^i-1)d]$ is right adjoint to $i_*$ by Theorem 2.1. Now we should only find an involution commuting with $i_*$. This time
the Kuhn duality does not work but we can take the ``Serre duality'' $D=(-)^{\#}\circ\Theta\circ\Theta$ instead. It commutes with ${\bf K}^r$
because the injective generators $S^{p^id}_U$ are self--dual. Hence we get the Poincar\'e duality.
Unfortunately the only easy to describe examples
are  $X=i^*(Y)$ for self--dual $Y\in {\cal DP}_d$, which explicitly give the Ext--groups
\[\mbox{Ext}^*_{{\cal P}_{p^id}}(P^{(i)},Y^{(i)}).\]
The problem is that it is difficult to find examples of self--dual $Y\in {\cal DP}_d$ other than injectives.\newline
It is also natural to try and find a version of Poincar\'e duality in which the simpler involution $(-)^{\#}\circ\Theta$ could play  role of $D$.
By Theorem 2.1 again, we conclude that if $p>2$ or $2|d$ then $D:=(-)^{\#}\circ\Theta$ commutes with ${\bf K}[\frac{(p^i-1)d}{2}]$.
With this formalism we can recover our Example 3.5, thus putting it into a wider context. We leave the details to the reader.

\end{document}